\documentclass[12pt]{article}
\usepackage{mathrsfs}
\usepackage{amsmath,amssymb}

\usepackage{amsfonts}
\usepackage{latexsym}
\usepackage{amsthm}
\usepackage{pictex}

\usepackage[T2A]{fontenc}     
\usepackage[cp1251]{inputenc} 
\usepackage{graphicx}

\def\lb{\label}

\newtheorem{theorem}{Theorem}[section]
\newtheorem{definition}{Definition}[section]
\newtheorem{lemma}{Lemma}[section]

\newtheorem{proposition}{Proposition}[section]

\begin{document}

\def\a{\alpha} \def\cA{{\cal A}} \def\bA{{\bf A}}  \def\mA{{\mathscr A}}
\def\b{\beta}  \def\cB{{\cal B}} \def\bB{{\bf B}}  \def\mB{{\mathscr B}}
\def\g{\gamma} \def\cC{{\cal C}} \def\bC{{\bf C}}  \def\mC{{\mathscr C}}
\def\G{\Gamma} \def\cD{{\cal D}} \def\bD{{\bf D}}  \def\mD{{\mathscr D}}
\def\d{\delta} \def\cE{{\cal E}} \def\bE{{\bf E}}  \def\mE{{\mathscr E}}
\def\D{\Delta} \def\cF{{\cal F}} \def\bF{{\bf F}}  \def\mF{{\mathscr F}}
\def\c{\chi}   \def\cG{{\cal G}} \def\bG{{\bf G}}  \def\mG{{\mathscr G}}
\def\z{\zeta}  \def\cH{{\cal H}} \def\bH{{\bf H}}  \def\mH{{\mathscr H}}
\def\e{\eta}   \def\cI{{\cal I}} \def\bI{{\bf I}}  \def\mI{{\mathscr I}}
\def\p{\psi}   \def\cJ{{\cal J}} \def\bJ{{\bf J}}  \def\mJ{{\mathscr J}}
\def\vT{\Theta}\def\cK{{\cal K}} \def\bK{{\bf K}}  \def\mK{{\mathscr K}}
\def\k{\kappa} \def\cL{{\cal L}} \def\bL{{\bf L}}  \def\mL{{\mathscr L}}
\def\l{\lambda}\def\cM{{\cal M}} \def\bM{{\bf M}}  \def\mM{{\mathscr M}}
\def\L{\Lambda}\def\cN{{\cal N}} \def\bN{{\bf N}}  \def\mN{{\mathscr N}}
\def\m{\mu}    \def\cO{{\cal O}} \def\bO{{\bf O}}  \def\mO{{\mathscr O}}
\def\n{\nu}    \def\cP{{\cal P}} \def\bP{{\bf P}}  \def\mP{{\mathscr P}}
\def\r{\rho}   \def\cQ{{\cal Q}} \def\bQ{{\bf Q}}  \def\mQ{{\mathscr Q}}
\def\s{\sigma} \def\cR{{\cal R}} \def\bR{{\bf R}}  \def\mR{{\mathscr R}}
\def\S{\Sigma} \def\cS{{\cal S}} \def\bS{{\bf S}}  \def\mS{{\mathscr S}}
\def\t{\tau}   \def\cT{{\cal T}} \def\bT{{\bf T}}  \def\mT{{\mathscr T}}
\def\f{\phi}   \def\cU{{\cal U}} \def\bU{{\bf U}}  \def\mU{{\mathscr U}}
\def\F{\Phi}   \def\cV{{\cal V}} \def\bV{{\bf V}}  \def\mV{{\mathscr V}}
\def\P{\Psi}   \def\cW{{\cal W}} \def\bW{{\bf W}}  \def\mW{{\mathscr W}}
\def\o{\omega} \def\cX{{\cal X}} \def\bX{{\bf X}}  \def\mX{{\mathscr X}}
\def\x{\xi}    \def\cY{{\cal Y}} \def\bY{{\bf Y}}  \def\mY{{\mathscr Y}}
\def\X{\Xi}    \def\cZ{{\cal Z}} \def\bZ{{\bf Z}}  \def\mZ{{\mathscr Z}}
\def\O{\Omega}
\def\ve{\varepsilon}
\def\vt{\vartheta}
\def\vp{\varphi}
\def\vk{\varkappa}

\newcommand{\gA}{\mathfrak{A}}
\newcommand{\gB}{\mathfrak{B}}
\newcommand{\gC}{\mathfrak{C}}
\newcommand{\gD}{\mathfrak{D}}
\newcommand{\gE}{\mathfrak{E}}
\newcommand{\gF}{\mathfrak{F}}
\newcommand{\gG}{\mathfrak{G}}
\newcommand{\gH}{\mathfrak{H}}
\newcommand{\gI}{\mathfrak{I}}
\newcommand{\gJ}{\mathfrak{J}}
\newcommand{\gK}{\mathfrak{K}}
\newcommand{\gL}{\mathfrak{L}}
\newcommand{\gM}{\mathfrak{M}}
\newcommand{\gN}{\mathfrak{N}}
\newcommand{\gO}{\mathfrak{O}}
\newcommand{\gP}{\mathfrak{P}}
\newcommand{\gR}{\mathfrak{R}}
\newcommand{\gS}{\mathfrak{S}}
\newcommand{\gT}{\mathfrak{T}}
\newcommand{\gU}{\mathfrak{U}}
\newcommand{\gV}{\mathfrak{V}}
\newcommand{\gW}{\mathfrak{W}}
\newcommand{\gX}{\mathfrak{X}}
\newcommand{\gY}{\mathfrak{Y}}
\newcommand{\gZ}{\mathfrak{Z}}

\def\mA{{\mathscr A}}
\def\mB{{\mathscr B}}
\def\mC{{\mathscr C}}
\def\mD{{\mathscr D}}
\def\mE{{\mathscr E}}
\def\mF{{\mathscr F}}
\def\mG{{\mathscr G}}
\def\mH{{\mathscr H}}
\def\mI{{\mathscr I}}
\def\mJ{{\mathscr J}}
\def\mK{{\mathscr K}}
\def\mL{{\mathscr L}}
\def\mM{{\mathscr M}}
\def\mN{{\mathscr N}}
\def\mO{{\mathscr O}}
\def\mP{{\mathscr P}}
\def\mQ{{\mathscr Q}}
\def\mR{{\mathscr R}}
\def\mS{{\mathscr S}}
\def\mT{{\mathscr T}}
\def\mU{{\mathscr U}}
\def\mV{{\mathscr V}}
\def\mW{{\mathscr W}}
\def\mX{{\mathscr X}}
\def\mY{{\mathscr Y}}
\def\mZ{{\mathscr Z}}

\def\Z{{\Bbb Z}}
\def\R{{\Bbb R}}
\def\C{{\Bbb C}}
\def\T{{\Bbb T}}
\def\N{{\Bbb N}}
\def\S{{\Bbb S}}
\def\H{{\Bbb H}}
\def\J{{\Bbb J}}
\def\dD{{\Bbb D}}

\def\qqq{\qquad}
\def\qq{\quad}
\newcommand{\ma}{\begin{pmatrix}}
\newcommand{\am}{\end{pmatrix}}
\newcommand{\ca}{\begin{cases}}
\newcommand{\ac}{\end{cases}}
\let\ge\geqslant
\let\le\leqslant
\let\geq\geqslant
\let\leq\leqslant
\def\ma{\left(\begin{array}{cc}}
\def\am{\end{array}\right)}
\def\iint{\int\!\!\!\int}
\def\lt{\biggl}
\def\rt{\biggr}
\let\geq\geqslant
\let\leq\leqslant
\def\[{\begin{equation}}
\def\]{\end{equation}}
\def\wh{\widehat}
\def\wt{\widetilde}
\def\pa{\partial}
\def\sm{\setminus}
\def\es{\emptyset}
\def\no{\noindent}
\def\ol{\overline}
\def\iy{\infty}
\def\ev{\equiv}
\def\/{\over}
\def\ts{\times}
\def\os{\oplus}
\def\ss{\subset}
\def\h{\hat}
\def\Re{\mathop{\rm Re}\nolimits}
\def\Im{\mathop{\rm Im}\nolimits}
\def\supp{\mathop{\rm supp}\nolimits}
\def\sign{\mathop{\rm sign}\nolimits}
\def\Ran{\mathop{\rm Ran}\nolimits}
\def\Ker{\mathop{\rm Ker}\nolimits}
\def\Tr{\mathop{\rm Tr}\nolimits}
\def\const{\mathop{\rm const}\nolimits}
\def\dist{\mathop{\rm dist}\nolimits}
\def\diag{\mathop{\rm diag}\nolimits}
\def\Wr{\mathop{\rm Wr}\nolimits}
\def\BBox{\hspace{1mm}\vrule height6pt width5.5pt depth0pt \hspace{6pt}}

\def\Diag{\mathop{\rm Diag}\nolimits}


\def\Twelve{
\font\Tenmsa=msam10 scaled 1200 \font\Sevenmsa=msam7 scaled 1200
\font\Fivemsa=msam5 scaled 1200 \textfont\msbfam=\Tenmsb
\scriptfont\msbfam=\Sevenmsb \scriptscriptfont\msbfam=\Fivemsb

\font\Teneufm=eufm10 scaled 1200 \font\Seveneufm=eufm7 scaled 1200
\font\Fiveeufm=eufm5 scaled 1200
\textfont\eufmfam=\Teneufm \scriptfont\eufmfam=\Seveneufm
\scriptscriptfont\eufmfam=\Fiveeufm}

\def\Ten{
\textfont\msafam=\tenmsa \scriptfont\msafam=\sevenmsa
\scriptscriptfont\msafam=\fivemsa

\textfont\msbfam=\tenmsb \scriptfont\msbfam=\sevenmsb
\scriptscriptfont\msbfam=\fivemsb

\textfont\eufmfam=\teneufm \scriptfont\eufmfam=\seveneufm
\scriptscriptfont\eufmfam=\fiveeufm}

\title {Schr{\"o}dinger operator on the zigzag
half-nanotube in magnetic field.  }

\author{
  Alexei Iantchenko
\begin{footnote}
{ Institute of Mathematics and Physics,
  Aberystwyth Univ.,
Penglais, Ceredigion,  SY23 3BZ, UK, email: aii@aber.ac.uk, {\em on leave from Malm{\"o} H{\"o}gskola, Sweden} }
\end{footnote} \and
Evgeny Korotyaev
\begin{footnote}
{School of Mathematics, Cardiff Univ., Senghennydd Road, Cardiff,
CF24 4AG, UK, e-mail: korotyaeve@cf.ac.uk}
\end{footnote}
}

\maketitle

\begin{abstract}
\no We consider the zigzag half-nanotubes
 (tight-binding approximation) in a uniform magnetic field which is described by the magnetic Schr\"odinger operator with a periodic potential plus  a finitely supported perturbation.
 We describe all eigenvalues and resonances of this operator, and
theirs dependence on the magnetic field.
  The proof is reduced to the analysis of    the periodic Jacobi
  operators on the half-line with finitely supported  perturbations.

\end{abstract}


\section {Introduction}
\setcounter{equation}{0}

After their discovery \cite{Ii}, carbon nanotubes remain in both theoretical and applied research (see \cite{SDD}). Structure of nanotubes are formed by rolling up a graphene sheet into a cylinder. Such nanomodels were introduced by Pauling \cite{Pa} in 1936 to simulate aromatic molecules. They were described in more detail by Ruedenberg and Scherr \cite{RS1} in 1953. Various physical properties of carbon nanotubes can be found in \cite{SDD}.

Single-wall nanotubes, one atomic layer in thickness in the radial direction,
are a very important variety of carbon nanotube because they exhibit important
electric properties that are not shared by the multi-walled carbon nanotube variants.
Single-wall nanotubes are the most likely candidate for miniaturizing electronics beyond
the micro electromechanical scale that is currently the basis of modern electronics.

 We consider the
Schr\"odinger operator $H^b=H_0^b+V+Q$  on the zigzag half-nanotube
$\G\ss\R^3$ (1D  tight-binding model of zigzag single-wall
half-nanotubes, see \cite{SDD}, \cite{N}) in a uniform magnetic
field $\mB=|\mB|{\bf e}_0,\ {\bf e}_0=(0,0,1)\in \R^3.$  Here
$H_0^b$ is the Hamiltonian of the nanotube in the magnetic field,
$V$ is the periodic potential of the nanotube, $Q$ is the finitely supported perturbation.

There are numerous mathematical results about Schr\"odinger operators on carbon nanotubes (zigzag, armchair and chiral) (see for example \cite{KL}, \cite{KL1}, \cite{K1}, \cite{KuP}, \cite{Pk}). All these papers consider the so called continuous models. But in the physical literature the most commonly used model is the tight-binding model.
In the tight binding model for a solid-state lattice of atoms, it is assumed that
the full Hamiltonian H of the system may be approximated by the Hamiltonian of an
isolated atom centered at each lattice point.
The mathematical models, e.g., the Schr\"odinger operator on the zigzag and
armchair nanotubes and ribbons in a uniform magnetic field
 $\mB$ and  in an external periodic electric potential were considered in \cite{KK1}, \cite{KK1}, \cite{Pk}, see also
\cite{RR}. For applications of our models see references in \cite{ARZ},
\cite{Ha}, \cite{SDD}.

 Our model
nanotube $\G$ is a graph (see Fig. \ref{f001} and 2) embedded in
$\R^3$ oriented in the $z$-direction ${\bf e}_0$ with unit edge
length. $\G$ is a set of vertices (atoms) ${\bf r}_{\o}$ connecting
by bonds (edges) $\G_{n,j,k}$ and
\begin{multline}
\G=\cup_{\o\in \cZ} {\bf r}_{\o},\qq {\bf r}_{n,0,k}={\bf
\vk}_{n+2k}+{3n\/2}{\bf e}_0,\ \qqq {\bf r}_{n,1,k}={\bf
r}_{n,0,k}+{\bf e}_0,\\
\omega=(n,j,k)\in\cZ=\Z_+\ts \{0,1\}\ts \Z_N, \qq \Z_N=\Z/(N\Z),
 \qq {\bf \vk}_k=R(\cos{\pi
k\/N},\sin{\pi k\/N},0),\\
\qq \ R={\sqrt 3\/4\sin {\pi\/2N}},\qqq \Z_+=\{j\in\Z,\,\,j\geq 0\}.
\end{multline}

\begin{figure}[h]
         \centering\includegraphics[clip]{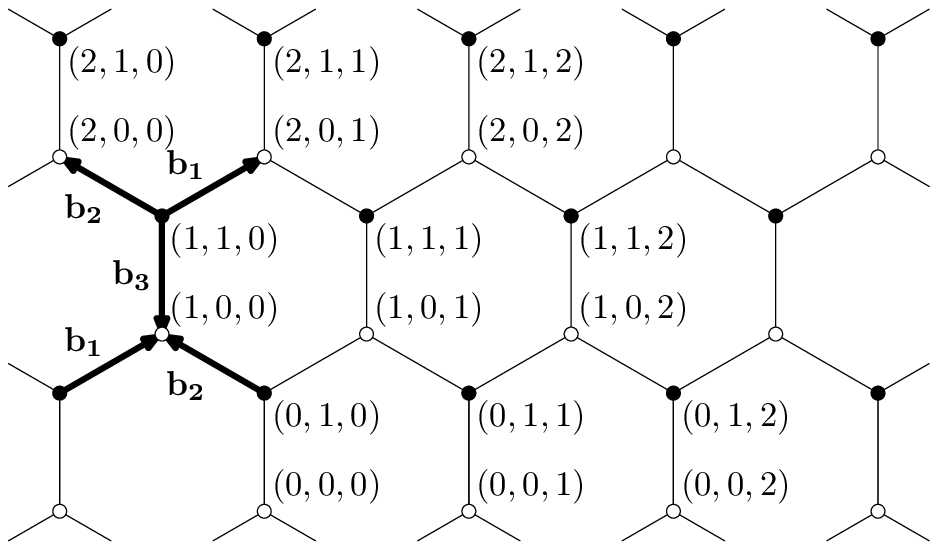}
         \begin{center}
         {\small Fig 1. A piece of a nanotube. }
         \end{center}
         \lb{f001}
\end{figure}
Our carbon model nanotube is the honeycomb lattice of a
graphene sheet rolled into a cylinder. This nanotube $\G$ has $N$
hexagons around the cylinder embedded in $\R^3$. Here $n \in \Z$
labels the position in the axial direction of the tube, $j=0,1$ is a
label for the two types of vertices (atoms)  (see Fig. \ref{f001}),
and  $k \in \Z_N$ labels the position around the cylinder. The
points ${\bf r}_{0,1,k}, k\in \Z_N$ are vertices of the regular
N-gon $\mP_0$ and ${\bf r}_{1,0,k}$ are the vertices of the regular
N-gon $\mP_1$. $\mP_1$ arises from $\mP_0$ by combination of the
rotation around the axis of the cylinder $\cC$ by the angle
${\pi\/N}$ and of the translation by  ${1\/2}{\bf e}_0$. Repeating
this procedure we obtain $\G$.

Introduce the Hilbert space $\ell^2(\G)$ of functions
$f=(f_\o)_{\o\in \cZ}$ on $\G$ equipped with the norm
$\|f\|_{\ell^2(\G)}^2=\sum_{\o\in \cZ} |f_\o|^2 $. The tight-binding
Hamiltonian $H^b$ on the half-nanotube $\G$ has the form
$H^b=H_0^b+\tilde{V}$ on $\ell^2(\G)$, where $H_0^b$ is  given by (see \cite{KL1})
\begin{align}
& (H_0^b f)_{n,0,k}=e^{ib_1}f_{n-1,1,k}+e^{ib_2}f_{n-1,1,k+1}+
 e^{ib_3}f_{n,1,k},\,\,\qq\,\,f_{-1,1,k}=0, \lb{010}\\
& (H_0^b f)_{n,1,k}=e^{ib_1}f_{n+1,0,k-1}+e^{ib_2}f_{n+1,0,k}+
 e^{-ib_3}f_{n,0,k},\qq f=(f_\o)_{\o\in \cZ},\nonumber\\
& \o=(n,j,k)\in \Z_+\ts \{0,1\}\ts \Z_{N},\qq
 b_3=0,\ \ b_1=-b_2=b={3|\mB|\/16}\cot {\pi\/2N},\nonumber
\end{align}
 and the
operator $\tilde{V}=V+Q$
 is given by
\begin{align}
& (\tilde{V}f)_\o=\tilde{V}_\o f_\o,\qq \mbox{where} \qq
\tilde{V}_{n-1,1,k}=\tilde{v}_{2n}, \qq
\tilde{V}_{n,0,k}=\tilde{v}_{2n+1},\qq \tilde{v}=(\tilde{v}_n)_{n\in\N}\in\ell^\iy,\lb{cb}\\
&\mbox{where}\,\,\tilde{v}_n=v_n+q_n\,\,\mbox{for}\,\,0\leq n\leq
p,\,\,q_p\neq 0,\,\,\mbox{and}\,\,\tilde{v}_n=v_n\,\,\mbox{for}\,\,
n>p.\nonumber
\end{align}

Such models can be realized using optical methods, by gating, or by an acoustic field (see
\cite{N}). For example, if an external potential is given by $A_0 \cos(\xi_0z + \beta_0)$ for some constant
$A_0,\xi_0,\beta_0,$ then we obtain
$$v_{2n}=A\cos\left(2\pi\xi \left(n-\frac13\right)+\beta\right),\qq v_{2n+1}=A\cos\left(2\pi\xi n+\beta\right),\qq n\in\N=1,2,\ldots,$$
for some constant $A,\xi,\beta.$ If $\xi$  is rational, then the sequence $v_n,$   is periodic for $n\in\N.$ If $\xi$ is
irrational, then the sequence $v_n$ $n\in Z_+$ is almost periodic.

We give the physical sense of the finitely supported potential $q=(q_n)_{n=0}^\iy,$ $q_p\neq 0$ and $q_n=0$
for all $n>p$. There are  two physical cases: a local defect in the nanotube and an effective potential.
The effective potential  is related to the boundary after cutting an infinite nanotube into  two pieces.
 The effective potential is due to an imperfection in the structure of the half-nanotube near the cut and corresponds to perturbations $q$ with  $p$ small. This motivates our detailed analysis of the properties of eigenvalues and resonances in the special case $p=1,$ $p=2,$ in Section \ref{S-examples}.

  In the present paper we suppose that the periodic background
  potential $v$ has period $2$ and is given by
  $v_{2n+1}=-v_{2n}=v\in\R,$ $n\in\N.$

\begin{figure}[h]
\lb{fig2}
         \centering\includegraphics[clip]{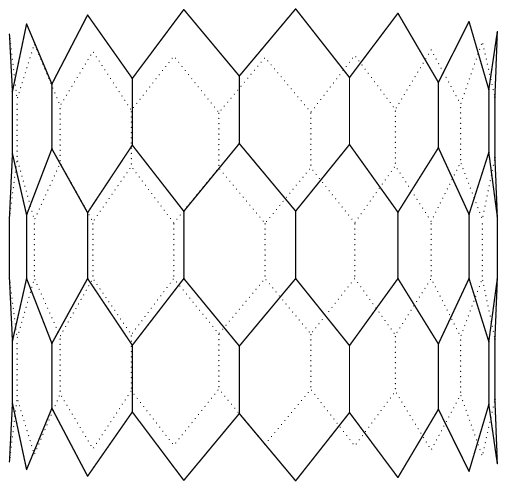}
         \begin{center}
         {\small Fig 2. Nanotube in the magnetic field. }
         \end{center}
\end{figure}

We formulate the result proven in \cite{KK2} in the form
convenient for us\\
 {\em
Each operator $H^b, b\in \R,$ is
unitarily equivalent to the operator $\os_1^N J_k^b$, where $J_k^b$
is the Jacobi operator, acting on $\ell^2(\N)$ and given by
\begin{align}
& (J_k^b y)_n=a_{n-1}y_{n-1}+a_{n}y_{n+1}+\tilde{v}_ny_n,\,\,
(\mbox{for}\,\,n\geq
 2),\qq (J_k^b y)_1=a_1y_2+\tilde{v}_1y_1\lb{Jk}\\
&
  a_{2n}\equiv a_{k,2n} =2|c_k(b)|,  \qq a_{2n+1}\equiv a_{k,2n+1}=1,\ \ c_k(b)=\cos (b+{\pi k\/N}), \
 n\in\N,\nonumber\\
 & \tilde{v}_n=v_n+q_n,\qq  q_j=0\,\,\mbox{for}\,\, j>p,\,\,
q_p\neq 0,
\end{align}
and $J_k^{b+{\pi\/N}}=J_{k+1}^b$, $J_k^{-b}=J_{N-k}^b$ for all
$(k,b)\in\Z_N\ts \R$. Moreover, the operators $H^{\pm b}$ and
$H^{b+{\pi\/N}}$ are unitarily equivalent for all $b\in \R$. }

{\no\bf Remarks.} 1) Note that the $n=1$ case in \ref{Jk} can be thought of as forcing the Dirichlet condition $y(0)=0.$ Thus, eigenfunctions must be non-vanishing at $n=1$ and eigenvalues must be simple.

 2)   The matrix of the operator $J_k^b$ is given by
\[
\lb{Jk1a}
 J_k^b=\left(\begin{array}{ccccc}
\tilde{v}_1   & 1     &0         & 0        &... \\
  1& \tilde{v}_2        &2|c_k|& 0        &... \\
  0  &2|c_k|& \tilde{v}_3         & 1     &... \\
  0  & 0        &1      & \tilde{v}_4       &... \\
  0  & 0        &0         &2|c_k|&... \\
  ... & ...      &...       &...       &... \\
 \end{array}\right).
\]

 If $c_k=\cos (b+{\pi k\/N})=0,$ then matrix (\ref{Jk1a}) has the form
\[\lb{spp0}
{J_k^b}|_{ c_k=0} =\cJ=\left(\begin{array}{ccccc}
 \tilde{v}_1   & 1     &0         & 0        &... \\
  1& \tilde{v}_2        &0& 0        &... \\
  0  &0& \tilde{v}_3         & 1     &... \\
  0  & 0        &1      & \tilde{v}_4       &... \\
  0  & 0        &0         &0&... \\
  ...& ...      &...       &...       &... \\
 \end{array}\right)=\os_{n\in \N}\cJ_n,\,\,\cJ_n=\left(
                                     \begin{array}{cc}
                                       \tilde{v}_{2n-1} & 1 \\
                                       1 & \tilde{v}_{2n} \\
                                     \end{array}
                                   \right),
\]
with the eigenvalues \[ \lb{spp}
\{z_{n,\pm}=v_n^+\pm |{v_n^-}^2+1|^{1\/2},\qq
v_n^\pm={\tilde{v}_{2n-1}\pm \tilde{v}_{2n}\/2} ,\qq n\in \N\}.
\]
Moreover, if $q_p\neq 0$, then\\
  if $p$ is
even,  then there are at most $p$ eigenvalues $z_{n,\pm},$ $n=1,2,\ldots, \frac{p}{2},$\\
 if $p$ is odd, then there are at most $p+1$  eigenvalues $z_{n,\pm},$ $n=1,2,\ldots,
\frac{(p-1)}{2}.$\\
Note that for some special choice of perturbations $\{q_1,\ldots q_p\}$ we can have $z_i^\pm=z_j^\pm$ for $i\neq j.$

 As perturbations have finite support, then  there are always two flat bands (two eigenvalues with infinite multiplicities) given by $z_{n,\pm}=\pm|{v}^2+1|^{1\/2},$ where $n\geq \frac{p}{2}+1$ if $p$ is even or  $n\geq
\frac{(p+1)}{2}+1$  if $p$ is odd.
 The flat bands are inherited from the pure periodic problem.

3) If $|c_k|={1\/2}$, then $J_k^b$ is the Schr\"odinger operator
with $a_{n}=1$ for all $n\in \N$. In particular, if $b=0,{N\/3}\in
\N$, then $J_{N\/3}^0$ is the Schr\"odinger operator.

4) Exner \cite{Ex} obtained some duality between Schr{\"o}dinger operators on graphs and certain Jacobi
matrices, which depend on energy. In our case the Jacobi matrices do not depend on energy.



{\bf  Unperturbed operator}. We start with the unperturbed operator
$H_0^b+V,$ which is unitary equivalent to $\os_1^N J_k^{b,0}.$ The
operator $J_k^{b,0}$ is acting in $\ell^2(\N)$ with Dirichlet boundary condition, see   (\ref{Jk}),
 where $\tilde{v}_n=v_n$  is the two-periodic
potential  verifying $v_{2n+1}=-v_{2n}=v\in \R.$
 It is known that, if $c_k\neq 0,$ then the absolutely
continuous spectrum of $J_k^{b,0}$ is given by two bands and the bound states in $\gamma_k^+$ (see
Section \ref{s-preliminaries}):
\begin{align}
&\s_{\rm ac}(J_k^{b,0})=[z_{k,0}^{b,+},z_{k,0}^{b,-}]\sm
\g_{k,1}^b,\qqq
\g_{k,1}^b= (z_{k,1}^{b,-},z_{k,1}^{b,+}),\label{abscont}\\
& z_{k,0}^{b,\mp}=\pm \sqrt{v^2+(2|c_k|+1)^2},\qqq
z_{k,1}^{b,\pm}=\pm \sqrt{v^2+(2|c_k|-1)^2},\qq k\in\Z_N, \nonumber\\
&\s(J_k^{b,0})=\s_{ac}(J_k^{b,0})\cup \s_{pp}(J_k^{b,0}), \qq
\s_{pp}(J_k^{b,0})=\left\{\begin{array}{lr}
                 \{v\}  & \mbox{if}\,\, \  1/2<|c_k|\leq 1, \\
                 \es & \mbox{if}\,\, \  0<|c_k|\leq 1/2,
                                \end{array}\right.\nonumber
\end{align}
where $\g_{k,1}^b$ is the middle gap in the spectrum of $J_k^{b,0}.$
We denote $\gamma_0=(-\infty,z_{k,0}^{b,+}),$
$\gamma_2=(z_{k,0}^{b,-},+\infty)$  the infinite gaps.

We denote $\Lambda=\L_k^b$ the two-sheeted Riemann surface for each
$J_k^b,$ obtained by joining the upper and low rims of two copies of
the cut plane $\C\setminus
 \s_{\rm
ac}(J_k^{b,0})$ in the usual (crosswise) way.   For $j=0,1,2,$ we
denote the copies of $\gamma_{j}$ on $\Lambda^+$ (respectively
$\Lambda^-$) by $\gamma_{j}^+$ (respectively $\gamma_{j}^-$), and
put  $\gamma_j^{\rm c}=\overline{\gamma_j^+}\cup
\overline{\gamma_j^-}.$ By abuse of notation we write also
$\gamma_j$ for $\gamma_j^+\cup \gamma_j^-$ and for its projection on
$\C.$


If $0<|c_k| < 1/2,$ then $v\in\gamma_1^-$ is an antibound state for
$J_k^b$ and if $|c_k|=1/2,$ then $v=z_{k,1}^{b,+}$ or
$v=z_{k,1}^{b,-}$ is virtual state (see Definition \ref{defstates}
below and Proposition \ref{Prop_states_nonperturbed}).

If  $c_k=0$ for some $(k,b)\in\Z_N\times\R,$ then (\ref{spp}) gives that the spectrum of $J_k^{b,0}$ is pure point:
$$\s (J_k^{b,0})=\s_{pp}(J_k^{b,0})=\{\pm\sqrt{v^2+1}\},$$ and each eigenvalue of $J_k^{b,0}$ is a flat band, i.e. has infinite multiplicity.

In \cite{KK2} it is shown that
 the spectral band $[z_{k,0}^{b,+}, z_{k,1}^{b,-}],$ (respectively $[z_{k,1}^{b,+}, z_{k,0}^{b,-}]$) shrinks to the flat  band $-\sqrt{v^2+1}$ (respectively $\sqrt{v^2+1}$) as $c_k\rightarrow  0$ and the
corresponding asymptotics are determined.


Let $b=\pi \left( \frac12 -\frac{1}{N}\right).$ Then $c_1=0$ and
$c_k=\cos\pi \left( \frac12 -\frac{1}{N}+\frac{
k}{N}\right)$   and the spectrum of $H^b$ is
given by
$$
\s(H^b)=\s_{ac}(H^b)\cup \s_{pp}(H^b), \qq \s_{pp}(H^b)=
    \{v,\pm \sqrt{1+v^2}\}  ,
$$
\[
\s_{ac}(H^b)=[z_{0}^{b,+},z_{0}^{b,-}]\sm \g(H^b),\ \g(H^b)=
(z_{1}^{b,-},z_{1}^{b,+}),
\]
where $\g(H^b)$ is the gap in the spectrum of $H^b$. If  $b\neq\pi \left( \frac12 -\frac{1}{N}\right)$  and all $c_k\neq 0,$
 $k=1,\ldots,N,$  then we obtain
$\{\pm\sqrt{1+v^2}\}\not\in\sigma_{pp}(H^b).$

Note that if $c_k=0$ for some $k\in\Z_N$ then $\s_{pp}(H^b)=\{\pm
\sqrt{1+v^2}\}\ss \g(H^b)$. From \cite{KK2} we know that
$\s(H^{b+{\pi\/N}})=\s(H^b )$ for all $b\in \R$. Then we need to
consider only the case $b\in [0,{\pi\/N})$ and in this case we get
$$
z_{0}^{b,+}=\ca z_{0,0}^{b,+}    & \mbox{if}  \qq b\le{\pi\/2N}\\
    z_{N-1,0}^{b,+}  & \mbox{if}  \qq b>{\pi\/2N}  \ac,
$$
Moreover, in particular case $\mB=0,$ ${N\/3}\in \N,\  b=0,$
 we obtain
$\g(H^0)=(-|v|, |v|).$

{\bf Finitely  supported perturbations.} We consider
the main  operator  $H^b=H_0^b+V+Q.$ Recall that $H^b$ is unitary equivalent to
$\os_1^N J_k^b,$ where $ J_k^b=J_k^{b,0}+q$ is given by
(\ref{Jk}) with $q_n=0$ for $n>p$ and the sequence
$y=(y_n)_{n=0}^\infty$ satisfies  the Dirichlet boundary condition
$y_0=0.$

  The perturbation $q$ does not change the  absolutely continuous spectrum:
 $\sigma_{\rm ac}(J_k^b)=\sigma_{\rm ac}(J_k^{b,0})=[\l_0^+,\l_1^-]\cup[\l_1^+,\l_0^-],$
where we used the simplified notations
 $\l_0^\pm\equiv z_{k,0}^{b,\pm}$ and
$\l_1^\pm\equiv z_{k,1}^{b,\pm}.$

 In our paper we  study the global properties of eigenvalues, virtual states and resonances of $J=J_k^b.$
Let $R(\l)=(J-\l)^{-1}$ denote the resolvent of  $J$ and  let
 $\langle\cdot,\cdot\rangle$ denote the scalar product in $\ell^2(\N).$
Then for any $f,g\in\ell^2(\N)$ the function $\langle Rf,g\rangle$ is defined on $\L_+$ outside the poles at the bound states $\l_0\in\gamma_j^+,$ $j=0,1,2.$ Recall that  the bound states are simple.
Moreover, if $f,g\in\ell^2_{\rm comp}(\N),$ where
  $\ell_{\rm comp}^2(\N)$ denotes the $\ell^2$ functions on $\N$ with finite support,
 then the function    $\langle Rf,g\rangle$  has an analytic extension from
  $\L_+$ into the Riemann  surface $\Lambda.$

 \begin{definition}\label{defstates} Let $c_k(b)\neq 0$ for some $b\in \R$.  \\
 1) A number $\l_0\in\Lambda_-$ is a resonance, if the
 function $\langle Rf,g\rangle$ has a pole at $\l_0$ for some $f,g\in\ell^2_{\rm comp}(\N).$ The
 multiplicity of the resonance is the multiplicity of the pole.
If $\Re \l_0 =0,$ we call $\l_0$  antibound state.\\
2) A real number $\l_0=\l_0^\pm$ or $\l_0=\l_1^\pm$ is a virtual state if
 $\langle Rf,g\rangle$ has a singularity at $\l_0$ for some $f,g\in\ell^2_{\rm comp}(\N).$\\
 3) The state  $\l\in\Lambda$  is a bound state or a resonance or a virtual state of $J.$

 We denote the set of all states of $J$ by $\gS\,(J).$

\end{definition}
In Section \ref{s-genJost}, \ref{defstates2},  we give an  equivalent characterization of the states.

 In the unperturbed case $J^0=J_k^{b,0}$ we
show in Proposition \ref{Prop_states_nonperturbed} that if
$0<c_k\leq 1,$ then $\gS\,(J^0)$ consists of one state: a bound state
$v\in\gamma_1^+,$ a antibound state $v\in\gamma_1^-$ or a virtual state
$v=\l_{1}^\pm$. Note that any such state is projected on the
Dirichlet eigenvalue $v\in\C,$ $\vp_2(v)=0.$
\\ \\  Let $\vartheta_n,\varphi_n$ be the fundamental solutions of the
  equation $ a_{n-1} y_{n-1}+a_{n}y_{n+1}+v_ny_n=\l y_n,$  satisfying
 $\vartheta_0=\varphi_1=1,$ $\vartheta_1=\varphi_0=0.$
Let  $f^\pm_n$ be the Jost
solution, $f^\pm_n=\tilde{\vartheta}_n+m_\pm\tilde{\varphi}_n,$ where
$\tilde{\vartheta}_n,$ $\tilde{\varphi}_n$  denote the solutions to (\ref{Jk}) satisfying
 $\tilde{\vartheta}_n=\vartheta_n,$ $\tilde{\varphi}_n=\varphi_n,$ for
 $n>p.$   Here $m_\pm$
  are the Titchmarch-Weyl functions. The functions $\vp,\vt$
  are polynomials, the Jost solutions $f^\pm$ and functions $m_\pm$ are
  meromorphic functions on $\Lambda.$ Note that $f^-(\l)=\overline{f^+(\overline{\l})},$ $\l\in\Lambda,$ and $f^\pm(\l)\in\ell^2(\N)$ for any $\l\in\L_\pm.$ We call $f_0^\pm$ the Jost functions.

We pass to the formulation of main results of the present paper. Recall that  all bound and virtual states of  $J\equiv J_k^b$ are simple (see Lemma \ref{L-virt}).
In the next theorem we give the characterization of the states of
$J_k^b.$
\begin{theorem}\lb{th-charact-states}
Let $c_k(b)\neq 0.$

i) The point $\l=v\in\gamma_1^+$ or $\l=v\in\gamma_1^-$  is a state of
$J=J_k^b$ iff  the projection of $\l$ on $\C$ is a zero of $\tilde{\varphi}_0.$
The value $\l\in\Lambda$ whose projection on the complex plane does
not coincide with $v$  is a state of $J$ iff $\l\in\Lambda$ is a
zero of the Jost function $f_0^+$:
$$\gS\,(J)\setminus \{v\}=\{\l\in\Lambda:\,\,f_0^+(\l)=0)\}\subset
\left(\cup_{j=0,1,2}\overline{\gamma_j^\pm}\right)\cup\Lambda_-.$$

ii) The state $\l=\l_{0,1}^\pm$ is a virtual state of $J$ iff one of the following two conditions is satisfied:\\
1) $\l\neq v$ and  $f_0^+(\l)=0;$\,\,  2)  $\l=v$ and
$\tilde{\varphi}_0(\l)=0.$

iii)   If $\l=v\in \gamma_1^-$ is an antibound
state for $J$ then it is necessarily simple.

 \end{theorem}





The distribution of the states  is summarized in
the following theorem.
\begin{theorem}\label{mainstatictheorem}
Let $c_k(b)\neq 0,$ $q_p\neq 0.$ Then the Jacobi operator $J\equiv J_k^b$ has $2p$ states
 counted with multiplicities. Moreover, the following facts hold true.\\
1) The total number of bound states and virtual states is  $\geq 2.$\\
2) In the closure of the middle gap $\gamma_1^{\rm c}=\overline{\gamma_1^+}\cup
\overline{\gamma_1^-}$ there is always an odd number of states  with at least one bound or virtual state.\\
3)  Let $\l_1<\l_2$ be any two bound
states of $J,$ $\l_{1,2}\in\gamma_k^+,$ for some  $k=0,1,2,$ such that there are no other
eigenvalues on the interval $\Omega^+=(\l_1,\l_2)\subset\gamma_k^+.$
Then there exists an odd number $\geq 1$ of
antibound states on $\Omega^-,$ where $\Omega^-\subset\gamma_k^-\subset\Lambda_-$ is the same interval
but on the second sheet.
\end{theorem}
{\bf Remarks.} 1)
If all $c_k\neq 0$ and $q_p\neq 0,$  then
the operator $\os_1^N J_k^b$ has in total $N 2p$ states.\\
2) If $p$ is even and
$q_1=q_3=\ldots=q_{p-1}=0,$ then $\{v\}$ is always a bound state or antibound state (see
Lemma \ref{Prop-l=v}).

In Theorem \ref{th-a0} we consider the limit of the states of each
$J_k^b$ as $c_k\rightarrow 0.$ Recall that  operator ${J_k^b}|_{ c_k=0}$  has two flat bands and a finite number of simple eigenvalues.
\begin{theorem}\label{th-a0} Let $z_{n,\pm},$ $n\in\N,$ be the eigenvalues  of the matrix ${J_k^b}|_{ c_k=0}$ given in
(\ref{spp0}). Let
 $c_k\rightarrow 0+.$  \\
 1) If $p$ is even, then\\
a) the set of  bound states of
$J_k^b$ converges to  the set ${\displaystyle
\{z_{n,\pm},\,\,n=1,\ldots, \frac{p}{2}}\}\ss\R,$\\
 b) the set of all resonances
of $J_k^b$ converges to the set of numbers\\
${\displaystyle\{z_{n,\pm},\,\, n=1,\ldots,
\frac{(p-2)}{2}}\}\cup\{\mu_{p-1}^0,\mu_p^0\},$ where only the numbers
\begin{equation}\lb{complex-resonances0}
\mu_{p-1,p}^0=v+\frac{q_{p-1}}{2}\pm\sqrt{\frac{q_{p-1}^2}{4}-\frac{q_{p-1}}{q_p}},\end{equation}
can be complex. \\
2) If $p$ is odd, then\\
a) the set of  bound states of $J_k^b$
converges to the set
${\displaystyle\{z_{n,\pm},\,\,n=1,\ldots, \frac{(p+1)}{2}}\}\ss \R;$\\
b) the
set of  resonances of $J_k^b$ converge to the set of real numbers
 ${\displaystyle\{z_{n,\pm},\,\, n=1,\ldots,
\frac{(p-1)}{2}} \}$.
\end{theorem}

In Theorem \ref{th_largeq} we consider the asymptotics of the states of
the half-nanotube Hamiltonian $H^b$ (unitary equivalent to $\os_1^N J_k^b$) for large
perturbation.
\begin{theorem}\label{th_largeq} Suppose  $q_j=q_j^0t,$ $j=1,\ldots,p,$   where all
$q_j^0\neq 0$  are fixed and $t>1.$ If $\l(t)\in\Lambda$ is a state of
$H^b,$ then either  $|\l(t)|\rightarrow \infty$
  or $\l(t)\rightarrow
(-1)^p v$  as $t\rightarrow\infty.$
\end{theorem}

If $v\rightarrow\infty,$ then (\ref{abscont}) implies that the absolutely continuous spectrum
degenerates  into two points $\{v\},$ $\{ -v\}.$

Suppose $p=2$ and $q_1=0$ and  $q_2$
is small enough, then the Hamiltonian  $H^b$
has
precisely $2N$ non-real complex conjugated resonances.
 More results about the cases  $p=1$ and $p=2$ are given in  Section
\ref{S-examples}.

The plan of the paper is as follows. In Section
\ref{s-preliminaries} we collect some well known facts about the
two-periodic Jacobi operators and its perturbations in the form
convenient for us.
\\
In Section \ref{s-genJost} we describe the properties of the perturbed operator.\\
In Section \ref{ssF} we consider the properties of the polynomial
$F=\varphi_2 f_0^+f_0^-$ which plays the crucial role in the proof
of the main results, similar to the case \cite{K2} .  Theorem \ref{mainstatictheorem} follows from
Lemma \ref{generalresults} and Theorem
 \ref{th-charact-states} follows from Lemmata \ref{L-virt} and \ref{L-ZerosF}. Theorems \ref{th-a0} and \ref{th_largeq}  follows from Lemmata  \ref{Jost-asymptotics} and \ref{l-asymptotics-F}. In Section
\ref{S-examples} we consider the cases $p=1$ and $p=2.$

\section{Periodic Jacobi operator.}\lb{s-preliminaries}
\setcounter{equation}{0}

In this section we recall some well known facts about  the infinite
Jacobi matrix $\J^0$
\[
\lb{21}
 \left(\begin{array}{cccccccc}
  ...& ...  & ...  & ...     &...         & ...        &... \\
 ...&a  & v  & 1     &0         & 0        &... \\
 ...&0  & 1& -v        &a& 0        &... \\
 ...&0  & 0  &a& v         & 1     &... \\
 ...&0  & 0  & 0        &1      & -v       &... \\
 ...&0  & 0  & 0        &0         &a&... \\
 ...&...  & ...& ...      &...       &...       &... \\
 \end{array}\right),\qq 0< a\leq 2
\]
and the associated equation for $ \J^0$
\[\lb{1e}
a_{n-1} y_{n-1}+a_{n}y_{n+1}+v_ny_n=\l y_n,\qq
 a_{2n+1}=1,\ a_{2n}=a,\,\,
 v_{2n+1}=v, v_{2n}=-v,
\]
$(\l,n)\in\C\ts\Z .$
 Introduce fundamental solutions
$\vp=(\vp_n(z))_{n\in\Z}$ and $ \vt=(\vt_n(z))_{n\in \Z} $ for
equation (\ref{1e}), under the condition $\vt_0=\vp_1=1$ and
$\vt_1=\vp_0=0$. We obtain
\begin{align}
& \vp_0=0, \qq \vp_1=1, \qq \vp_2={\l}-v,\qqq \vp_3={\l^2-v^2-1\/a},
 \nonumber\\
& \vt_0=1, \qq \vt_1=0, \qq \vt_2=-a, \qq \vt_3=-\lambda -v \qq
..... \label{Weyl-solutions-short}
\end{align}
 The monodromy matrix $M_2$
satisfies
\[
\lb{M2}
M_2(\l)=  \ma \!\!\!\vt_{2} & \vp_{2} \!\!\!\\
           \!\!\!\vt_{3} & \vp_{3}\!\!\!
       \am =\ma -a & \l-v \\ -\l-v & {\l^2-v^2-1\/a} \am .
\]
The Lyapunov function is defined in the standard way:
\[
\lb{D2} \D={\Tr M_2\/2}={\l^2-v^2-a^2-1\/2a}=\cos 2\vk ,
\] where $\vk$ is the Bloch quasimomentum.

The periodic eigenvalues $\l_{0}^{\pm}$ satisfy the equation
$\D(\l)=1$ and the anti-periodic eigenvalues $\l_{1}^{\pm}$ satisfy
the equation $\D(\l)=-1$ and they are given by
\[
\lb{eib}
 \l_{0}^{\mp}=\pm \sqrt{v^2+(a+1)^2},\qqq
\l_{1}^{\pm}=\pm \sqrt{v^2+(a-1)^2}.
\]

 The absolutely continuous spectrum of $\J^0$ has the form
\[
\lb{eib2} \s_{\rm ac}(\J^0)=[\l_{0}^{+},\l_{1}^{-}]\cup
[\l_{1}^{+},\l_{0}^{-}]=[\l_{0}^{+},\l_{0}^{-}]\sm \g_{1},\qq
\g_{1}=(\l_{1}^{-},\l_{1}^{+})
\]
where $\g_{1}$ is a gap.  Note that
$
  \g_{1}=(\l_{1}^{-},\l_{1}^{+})\ne \es,\qq \mbox{if} \qq
|v|+|a-1|>0.
$
We denote also $\gamma_0=(-\infty,\lambda_0^+)$ and
$\gamma_2=(\lambda_0^-,+\infty).$



We recall from the Introduction that the two-sheeted Riemann surface $\Lambda$  is obtained by joining the
upper and low rims of two copies  $\Lambda^\pm$ of the cut plane
$\C\setminus
 \s_{\rm
ac}(\J^0)$ in the usual (crosswise) way.   For $j=0,1,2,$  $\gamma_{j}^+$ (respectively $\gamma_{j}^-$) denote
the copies of $\gamma_{j}$ on $\Lambda^+$ (respectively
$\Lambda^-$), and
 $\gamma_j^{\rm c}=\overline{\gamma_j^+}\cup
\overline{\gamma_j^-}.$ By abuse of notation we write also
$\gamma_j$ for $\gamma_j^+\cup \gamma_j^-$ and for its projection on
$\C.$

The eigenvalues of $M_2$ are given by
$\xi_\pm^2=\Delta\pm\sqrt{\Delta^2-1}.$ On $\gamma_0^+,$ we choose
\[\lb{multipliers}\mbox{for}\,\, \lambda\in \gamma_0^+=(-\infty, \lambda_0^+)\subset\Lambda^+,\,\, \xi_+^2=
\Delta-\sqrt{\Delta^2-1},\,\,\xi_-^2=\Delta+\sqrt{\Delta^2-1}.
\]
 For
others $\l\in\Lambda,$ the functions $\xi_\pm(\lambda)$ are defined
by an analytic continuation.

If $\lambda=\pm v\in\gamma_1^+$  (these numbers will play a special
role later)  then $\Delta (\pm v)=\frac{-a^2-1}{2a}$ and
\begin{equation}\label{xi+v}\xi^2_+(\pm v)=\frac{-a^2-1}{2a}
+\left|\frac{a^2-1}{2a}\right|=\left\{\begin{array}{ll}
                                                                    -a & \mbox{if}\,\, 0<a<1, \\
                                                                    -1/a & \mbox{if}\,\,
                                                                    a>1,
                                                                  \end{array}\right.\end{equation}
                                                                  and
                                                                  opposite
                                                                  for
                                                                  $\xi_-^2(\pm
                                                                  v).$

 Then for $\l\in\gamma_0\cup\gamma_1\cup\gamma_2$ we have $|\xi_+^2| <1$ and $|\xi_-^2| >1.$
The eigenvectors of $M_2$ are chosen in the form $(1,m_\pm)$ and then
the Titchmarsh-Weyl functions are \[\lb{Weyl-function}
m_\pm(\l)=\frac{\xi_\pm^2-\vartheta_2}{\varphi_2}=\frac{\xi_\pm^2+a}{\lambda
-v}.\]
 For
$\lambda\in\gamma_1^+$ we have also
\begin{align}
m_\pm&=
\frac{\phi\pm\sqrt{\Delta^2-1}}{\varphi_2}=\frac{\Delta+a\pm\sqrt{\Delta^2(\lambda)-1}}{\l
-v}={\f\pm i\sin 2\vk\/\l-v}\lb{m}\\
 \f
 &={\vp_3-\vt_2\/2}={\l^2-v^2+a^2-1\/2a}=\Delta+a.\label{phi-function}
\end{align}

 On each $\gamma_k^+,$ $k=0,1,2,$ the quasimomentum $\varkappa(\l)$
 has constant positive imaginary part and we put $\varkappa=i h,$
 $h=h_k>0.$ Then $\Delta=\cosh (2h)$ and
 \[\lb{isin}
i\sin 2\varkappa=-(-1)^{k}\sqrt{\Delta^2(\l)-1}=-(-1)^{k}\sinh 2h.\]

Now the Floquet solutions $\p_n^\pm =\vt_n+m_\pm \vp_n$ are
\[\label{Floquet-solutions}
\p_0^\pm=1,\qq \p_1^\pm=m_\pm, \qq \p_2^\pm=e^{\pm
2i\vk}=\x_\pm^2,\qqq \p_{2n}^\pm=\x_\pm^{2n}, \qqq
\p_{2n+1}^\pm=\x^{2n}_\pm m_\pm,
\]
where $\xi_\pm^2=e^{\pm 2 i \vk}$ are the Floquet multipliers. Recall that $\psi_n^\pm\in\ell^2(\N)$ for any $\l\in\Lambda^\pm.$

Note the following simple identities which will be used in the paper:
\begin{equation}\label{Fact7}
\phi^2+1-\Delta^2=1-\varphi_3\vartheta_2=-\vartheta_3\varphi_2.
\end{equation}
Let $\{\phi_n,\psi_n\}=a_n(\phi_n\psi_{n+1}-\phi_{n+1}\psi_n\}$
denote the  Wronskian.

In the next theorem we describe  the states of the restriction of $\J^0$ to $\N$ defined in (\ref{Jk}) with $\tilde{v}_n=v_n$.

\begin{proposition}[Unperturbed case]\lb{Prop_states_nonperturbed}
The half-periodic Jacobi operator $J^0$   given by equation (\ref{Jk})  with $\tilde{v}_n=v_n,$
 has absolutely continuous spectrum (\ref{eib2}):
$\sigma_{\rm ac}\,(J^0)=[\l_{0}^{+},\l_{1}^{-}]\cup
[\l_{1}^{+},\l_{0}^{-}]$ and a state at $\l =v\in\overline{\gamma_1^+}\cup \overline{\gamma_1^-},$ whose
projection $v\in\C$ satisfies $\varphi_2(v)=0.$ There are three
possibilities:\\
 if $a>1$ then $J^0$ has  simple bound state at $\l=v\in\gamma_1^+;$\\
if $0<a<1$ then $J^0$ has simple antibound state at $\l=v\in\gamma_1^-;$\\
if $a=1$ then $\l=v$ is a simple virtual state, $v=\l_1^+$ or
$v=\l_1^-$ if $v>0$ respectively $v<0.$
\end{proposition}

{\bf Proof:}
 The kernel of the resolvent of $J^0$  is given by
$$R_0(n,m)=-\frac{\varphi_n\psi_m^+}{\left\{\varphi,\psi^+\right\}},\,\,n
<m,$$ where  $\left\{\varphi,\psi^+\right\}=-a.$ According to Lemma \ref{defstates2} (see Section \ref{s-genJost}), the bound states
(resonances) are the poles of
${\mR}_0(n)=\psi^+_n(\lambda)=\vartheta_n(\lambda)+m_+(\l)\varphi_n(\l)$
 on $\Lambda_+$ (respectively on $\Lambda_-$).
Hence, the only state is the pole of $m_+$  on $\Lambda_\pm,$ whose
projection on $\C$ is the zero of $\varphi_2(\lambda),$  i.e. $\l
=v\in\gamma_1.$

We have
$$m_+=\frac{\xi_+^2+a}{\l -v},\,\,a=2|c_k|,\,\,c_k=\cos\left(
b+\frac{\pi k}{N}\right).$$

If $0<a<1,$ then by  (\ref{xi+v})) $\l =v\in\gamma_1^+$ is a simple
zero for the numerator while at $\l =v\in\gamma_1^-$ the numerator
is non-zero. Thus $\l=v$ is an antibound state. Similar we get that if
$1<a<2,$  then $\l=v$ is a bound state.

If $a=1$ then $\Delta=(\l^2-v^2-2)/2$ and
$$\Delta^2-1=-(\l -v)(\l+v) +\frac{(\l-v)^2(\l +v)^2}{4}.$$
Suppose $v>0,$ then $v=\l_1^+.$ Let $\l-v=-\epsilon,$ $\epsilon >0,$
and let $\epsilon\rightarrow 0.$ Then
$$\Delta=-1-v\epsilon+{\mathcal O}(\epsilon^2),\,\,\sqrt{\Delta^2-1}=\sqrt{\epsilon}\sqrt{2v}+{\mathcal
O}(\epsilon),$$ and
\begin{equation}\lb{unpertvirt}
m_+(v-\epsilon)=\frac{\Delta+a+\sqrt{\Delta^2-1}}{\l
-v}=\frac{\sqrt{\epsilon}\sqrt{2v} +{\mathcal
O}(\epsilon)}{\epsilon}=\frac{\sqrt{2v}}{\sqrt{\epsilon}} +{\mathcal
O}(1).\end{equation} Thus if $a=1,$ the function $(\mR_n(\l))^2$ has
a pole at $\l=v$ and $\l=v$ is a virtual state. \qed

\section{ Jost
functions}\lb{s-genJost}
\setcounter{equation}{0}


 We introduce  the  Jost solutions
 as  solutions $f_n^\pm,$ of the equation
  \begin{equation}\lb{21.5}
 a_{n-1}
y_{n-1}+a_{n}y_{n+1}+\tilde{v}_ny_n=\l y_n,\qq
n\in\N,\qq\l\in\Lambda,
\end{equation} satisfying
\begin{equation}\label{defJost}
f_n^\pm=\p_n^\pm,\,\,  \mbox{for}\,\, n > p,
\end{equation}
where $\psi_n^\pm$ are the Floquet solutions (\ref{Floquet-solutions}) for the unperturbed problem,  and  $\tilde{v}_j=v_j+q_j$  with $q_n=0$
for $n>p.$
We recall that, as in (\ref{1e}), we have
  $v_{2n+1}=-v_{2n}=v\in \R,$ $
    a_{2n+1}=1,\ a_{2n}=a=2|c_k|\ne 0,\,\,
    c_k=\cos (b+{\pi k\/N}).$
 We have
$\overline{f^\pm_n}(\overline{\lambda})=f^\mp_n(\lambda),$ $\l\in\Lambda.$

The equation (\ref{21.5}) has unique solutions
$\tilde{\vartheta}_n,$ $\tilde{\varphi}_n$ such that
$$\tilde{\vartheta}_n(\lambda)=\vartheta_n(\lambda),\,\,\tilde{\varphi}_n(\lambda)=\varphi_n(\lambda)\,\,\mbox{for}\,\,n>p,\,\,\lambda\in
\C.$$ The functions $\tilde{\vartheta}_n(\cdot),$
$\tilde{\varphi}_n(\cdot)$ are polynomials. The functions $f_n^\pm$
have the form
\begin{equation}\label{3.2}
f^\pm_n=\tilde{\vartheta}_n +m_\pm \tilde{\varphi}_n,\,\,m_\pm={\f
\pm i\sin
2\vk\/\vp_2}=\frac{\Delta+a\pm\sqrt{\Delta^2(\lambda)-1}}{\l -v}.
\end{equation}
 Here $\f$ is defined in (\ref{phi-function}),
 $\varphi_2=\l -v$ and $\Delta$ is the Lyapunov function. The
functions $f^\pm_0$ are called Jost functions. The Jost functions
are analytic at all $\l\in\Lambda$ whose projection on the complex
plane $\C$ is different from $v,$   and has branch points
$\l_{0,1}^\pm.$

The asymptotics of the Jost functions are given in the following
Lemma.

\begin{lemma}\lb{Jost-asymptotics} Let $p,n\in\N$ and $p\geq n.$  Suppose $q_p\neq 0.$\\
{\bf 1)} If $p$ is even ($a_p=a,\,\,v_p=-v$), then\\ \\for
$\l\in\gamma_{0,2}^+$ in the limit $|\l|\rightarrow\infty,$ we have
\begin{align*}
 f_{0}^+&=1-\lambda^{-1}\sum_{k=1}^pq_k
+{\mathcal O}(\lambda^{-2}),\qqq
 f_{0}^-=
\frac{\l^{2p-1}}{a^p} \left[-q_p +{\mathcal
O}(\lambda^{-1})\right];
\end{align*}
for $\l\in\gamma_1^+$ as $a\rightarrow 0+,$ we have
\begin{align}
f_0^+&=(2\delta)^{-p/2}
\prod_{k=1}^{p-1}\rt[(\l-\tilde{v}_{k})(\l-
\tilde{v}_{k+1})-1\rt]+\delta^{-p/2}{\mathcal O}(a^2) ,\nonumber\\
f_0^-&=\frac{(2\delta)^{p/2}}{a^p}\left[(\l-\tilde{v}_{p-1})\left\{(\l-\tilde{v}_p)-
\frac{2\delta}{\l-v}\right\}-1\right]\prod_{k=1}^{p-3}
\rt[(\l-\tilde{v}_{k})(\l-
\tilde{v}_{k+1})-1\rt]\lb{complex-resonances1}\\&\hspace{13cm}+\delta^{p/2}{\mathcal
O}(a^2) ,\nonumber
\end{align}
where $\delta=(\l^2-v^2-1)/2;$\\ if  $q_k=tq_k^0$ with all $q_k^0\neq 0$
and  $t\rightarrow\infty,$ then we have
 ${\displaystyle f_0^+=t^p{\x^p_+\/a^{p/2}}\prod_{k=1}^{p}q_k^0+{\mathcal O}(t^{p-1}).}$\\ \\
{\bf 2)} If $p$ is odd ($a_p=1,$ $v_p=v$), then \\ \\
for $\l\in\gamma_{0,2}^+$ in the limit $|\l|\rightarrow\infty,$ we have
\begin{align}
f_0^+&=1-\lambda^{-1}\sum_{k=1}^pq_k+{\mathcal O}(\lambda^{-2}),\qqq
 f_0^-=
\frac{\l^{2p-1}}{a^{p+1}}\left[ -q_p+ {\mathcal
O}(\l^{-1})\right];
\end{align}
for  $\l\in\gamma_1^+,$ as $a\rightarrow 0+,$ we have
\begin{align*}
f_{0}^+&=(2\delta)^{-(p+1)/2}\left[(\l -\tilde{v}_p)(\l+v) -1\right]\prod_{k=1}^{p-2}\left[(\l-\tilde{v}_{k})(\l-\tilde{v}_{k+1})-
1\right]+\delta^{-(p+1)/2}{\mathcal
O}(a^2),\\
f_{0}^-&=\frac{(2\delta)^{(p+1)/2}}{a^{p+1}}\cdot\frac{ -q_p}{\l-v}\prod_{k=1}^{p-2}\left[(\l-\tilde{v}_{k})
(\l-\tilde{v}_{k+1})-1\right] +\delta^{(p+1)/2}{\mathcal
O}(a^2),
\end{align*}
if  $q_k=tq_k^0$ with $q_k^0\neq 0$ and  $t\rightarrow\infty,$ then
 ${\displaystyle f_0^+=t^p{\x^p_+\/a^{(p+1)/2}}\frac{1+a\xi_-^2}{\l-v}\prod_{k=1}^pq_k^0 +
 {\mathcal O}(t^{p-1}).}$
\end{lemma}

The proof  is technical  and uses the standard arguments.
The asymptotics of $f_0^+$ on $\gamma_{0,2}^+$ as  $\l\rightarrow\infty$ are well known (see  for example Teschl
 \cite{T}).

It is well known that the spectrum of $J=J_k^b,$ introduced in (\ref{Jk}), consists of absolutely
continuous part $\sigma_{\rm ac}\,(J)=\sigma_{\rm ac}\,(J^0)$ and a
finite number of simple bound states in each gap $\gamma_{k}^+,$
$k=0,1,2.$ The states of $J$  correspond to the poles of a
meromorphic function: resolvent or its square.

  The kernel of the resolvent of $J$ is
$$R(n,m)=\langle e_n,(J-\lambda)^{-1} e_m\rangle=-\frac{\Phi_n f_m^+}{\left\{\Phi,f^+\right\}},\,\,n
<m,$$ where $e_n=(\delta_{n,j})_{j\in\N},$ $J\Phi_n=\l\Phi_n,$
$\Phi_0=0,$ $\Phi_1=1,$ and the Wronskian $\left\{\Phi,f^+\right\}=-a_0f^+_0.$

Each function $\Phi_n(\lambda),$ $n\in \N,$ is polynomial in
$\lambda.$
The function $R(n,m)$ is meromorphic on $\Lambda$ for each
$n,m\in\Z.$ The singularities of $R(n,m)$ are given by the
singularities of
$${\mR}_n(\l)=\frac{f^+_n(\lambda)}{f^+_0(\lambda)}=
\frac{\tilde{\vartheta}_n(\lambda)+m_+(\l)\tilde{\varphi}_n(\l)}{f^+_0(\lambda)}.$$
The following Lemma follows from  Definition
\ref{defstates}.
\begin{lemma}\label{defstates2}
 1) A real number $\l_0\in\gamma_j^+,$ $j=0,1,2,$ is a bound state, if the
 function ${\mR}_n(\l)$ has a pole at $\l_0$ for almost all
 $n\in\N$ (eventually except a finite number of $n$'s) (it is known that the bound states are simple). \\
 2) A number $\l_0\in\Lambda_-,$ is a resonance, if the
 function ${\mR}_n(\l)$ has a pole at $\l_0$ for almost all
 $n\in\N$ (eventually except a finite number of $n$'s). The
 multiplicity of the resonance is the multiplicity of the pole.
If $\Re \l_0 =0,$ we call $\l_0$  antibound state.\\
3) A real number $\l_0=\l_0^\pm$ or $\l_0=\l_1^\pm$ is a virtual state if
$({\mR}_n(\l))^2$ or ${\mR}_n(\l)$ has a pole at $\l_0$ for almost
all
 $n\in\N$ (eventually except a finite number of $n$'s).\\
 4)  The state  $\l\in\Lambda$  is a bound state, resonance or virtual state.\\

\end{lemma}
We recall that  the set of all states of $J$ is denoted by $\gS\,(J).$

 Each function $f^+_n(\lambda),$ $n\in\N,$ is
analytic at all $\l\in\Lambda$ whose projection on the complex
plane $\C$ is different from $v.$  The Jost function $f^+_0(\l)$ has
finite number of real zeros on each $\gamma_k^\pm$ and finite number
of complex conjugated zeros on $\Lambda_-.$

Remark that  if $\l_0\in\gamma_k,$ for some $k=0,1,2$
(then $\l_0\neq\l_{0,1}^\pm$), and if $f_0^+(\l_0)\neq 0$ and
$\varphi_2(\l_0)\neq 0,$ then the resolvent is  analytic at $\l_0.$

 As
$\tilde{\vartheta}_n(\lambda),$ $\tilde{\varphi}_n(\l)$ are
polynomials then the singularities are zeros of $f_0^+(\l)$ and
eventually  singularities of $m_+$ at $\l=v$ (as in the unperturbed
case).

To describe the states of the general operator $J$ stated in Theorem
\ref{mainstatictheorem}, it is convenient to introduce a special
polynomial
 whose zeros give all  states of $J.$

\section{Function $F$ and  proofs of main results}\label{ssF}
\setcounter{equation}{0}
We introduce function
$F(\lambda)=\varphi_2f^+_0f^-_0.$
\begin{lemma}\lb{l-asymptotics-F}Suppose $q_p\neq 0.$\\  i) The function $F(\lambda)=\varphi_2f^+_0f^-_0(\l)$
is polynomial of degree $ 2p$ and satisfies
\begin{equation}\lb{formulaF}
F=\varphi_2\tilde{\vartheta}_0^2+2\phi\tilde{\vartheta}_0\tilde{\varphi}_0-\vartheta_3\tilde{\varphi}_0^2=(\l
-v)\tilde{\vartheta}_0^2+\frac{1}{a}(\l^2-v^2+a^2
-1)\tilde{\vartheta}_0\tilde{\varphi}_0+(\l +v)\tilde{\varphi}_0^2.
\end{equation}
ii) For $\l\in\C,$ in the limit
$|\l|\rightarrow\infty,$ we have asymptotics
\begin{align}
F &= \frac{\l^{2p}}{(a_p\ldots a_{0})^2} \left[-a^2q_p +{\mathcal
O}(|\lambda|^{-1})\right],\,\,\,\,\mbox{if}\,\,p\,\,\mbox{is even},\label{Feven-large-lambda}\\
F &= \frac{\l^{2p}}{(a_p\ldots a_{0})^2} \left[-q_p +{\mathcal
O}(|\lambda|^{-1})\right],\,\,\,\,\mbox{if}\,\,p\,\,\mbox{is
odd},\label{Fodd-large-lambda}
\end{align}
where ${\mathcal O}(|\lambda |^{-1})$ is uniformly bounded in $a.$
In particular,  if $\l\in\R$ and $|\l|\rightarrow\infty,$ we have
$\sign(F)=-\sign(q_p).$

iii) In the limit $a\rightarrow 0,$ the function $F$ behaves as
follows:\\
if $p$ is even, then
\begin{align}
F=&\frac{1}{a^p}\left[(\l-\tilde{v}_{p-1})(\l-\tilde{v}_p)-
1 \right]\left[(\l-\tilde{v}_{p-1})(1-q_p(\l-v))-(\l-v) \right]\lb{complex-resonances2}\\ &\cdot\prod_{k=1}^{p-3}\left[(\l-\tilde{v}_{k})(\l-
\tilde{v}_{k+1})-1\right]^2
+{\mathcal O}(a^{2-p}),
\end{align}
if $p$ is odd
\begin{align*}
&F=\frac{-q_p}{a^{p+1}}\left[(\l-\tilde{v}_{p})(\l+v)- 1
\right]
\prod_{k=1}^{p-2}\left[(\l-\tilde{v}_{k})(\l-
\tilde{v}_{k+1})-1\right]^2 +{\mathcal O}(a^{1-p}).
\end{align*}
Here  ${\mathcal O}(a^j)$ is uniformly bounded in $\lambda\in\C.$

iv) Put $q_j=tq_j^0$ for all $q_j^0\neq 0$ fixed and
$t\rightarrow\infty.$ Then,
\begin{align*}
&F(\l)=t^{2p}\frac{(\l-v)}{a^p}\left[\prod_{k=1}^p(q_p^0)^2
+{\mathcal O}(t^{2p-1})\right],\,\,\,\,\mbox{if $p$ is even}\\
&F(\l)=t^{2p}\frac{(\l+v)}{a^{p+1}}\left[\prod_{k=1}^p(q_p^0)^2
+{\mathcal O}(t^{2p-1})\right],\,\,\,\,\mbox{if $p$ is odd},
\end{align*}
uniformly bounded in $\l\in\C.$

\end{lemma}
{\bf Proof:} We have
\begin{align*}f^+_0f^-_0&=(\tilde{\vartheta}_0+m_+\tilde{\varphi}_0)(\tilde{\vartheta}_0+m_-\tilde{\varphi}_0)=\tilde{\vartheta}_0^2+(m_++m_-)\tilde{\vartheta}_0\tilde{\varphi}_0+m_+m_-\tilde{\varphi}_0^2\\
&=\tilde{\vartheta}_0^2+\frac{2\phi}{\varphi_2}\tilde{\vartheta}_0\tilde{\varphi}_0+\frac{\phi^2+1-\Delta^2}{\varphi_2^2}\tilde{\varphi}_0^2=
\tilde{\vartheta}_0^2+\frac{2\phi}{\varphi_2}\tilde{\vartheta}_0\tilde{\varphi}_0+\frac{-\vartheta_3\varphi_2}{\varphi_2^2}\tilde{\varphi}_0^2\\
&=
\tilde{\vartheta}_0^2+\frac{2\phi}{\varphi_2}\tilde{\vartheta}_0\tilde{\varphi}_0-\frac{\vartheta_3}{\varphi_2}\tilde{\varphi}_0^2,
\end{align*}
where we have used (\ref{m}) and (\ref{Fact7}). The degree $2p$ will
come from  $\varphi_2=\lambda -v$ and asymptotics
(\ref{Feven-large-lambda}), (\ref{Fodd-large-lambda}).

Now as $F$ is polynomial, in order to prove the asymptotics
$|\l|\rightarrow \infty$ on $\C$ it is enough to consider
$\l\rightarrow +\infty$ on $\gamma_2^+,$ $a\rightarrow 0$ or
$q_j\rightarrow\infty$ for $\l\in\gamma_1^+.$ The proof thus follows
from the asymptotics of the  Jost functions given in Lemma
\ref{Jost-asymptotics}.  \qed

From iii), Lemma \ref{l-asymptotics-F}, we get the leading orders of the zeros of $F$ as $a\rightarrow 0$ which correspond to the leading orders of the states.   Using Lemma \ref{Jost-asymptotics} we know if the limiting state is a bound state or a resonance.  Recall the eigenvalues of the matrix (\ref{spp0}) given by $z_{n,\pm}=v_n^+\pm|{v_n^-}^2+1|^{1\/2}$  (see (\ref{spp})). If $q_p\neq 0$ and\\
  if $p$ is
even,  then there are at most $p$ eigenvalues $z_{n,\pm},$ $n=1,2,\ldots, \frac{p}{2},$ where
$$
v_n^+={\tilde{v}_{2n-1}+ \tilde{v}_{2n}\/2}={q_{2n-1}+q_{2n}\/2},\qqq v_n^-={\tilde{v}_{2n-1} - \tilde{v}_{2n}\/2}=v+{q_{2n-1} -
q_{2n}\/2}
,$$\\
 if $p$ is odd, then there are at most $p+1$  eigenvalues $z_{n,\pm},$ $n=1,2,\ldots,
\frac{(p-1)}{2},$ where
$$ v_n^+={\tilde{v}_{2n-1}+ \tilde{v}_{2n}\/2}={q_{2n-1}+q_{2n}\/2}
,\,\,v_n^-={\tilde{v}_{2n-1} - \tilde{v}_{2n}\/2}=v+{q_{2n-1} - q_{2n}\/2},$$
$$ v_{(p+1)/2}^+={\tilde{v}_p-v\/2}={q_p\/2}
,\,\,v_{(p+1)/2}^-={\tilde{v}_p+v\/2}=v+\frac{q_p}{2}.$$
 Recall that, as perturbations have finite support, then  there are also two flat bands (two eigenvalues with infinite multiplicities) given by $z_{n,\pm}=\pm|{v}^2+1|^{1\/2},$ where $n\geq \frac{p}{2}+1$ if $p$ is even or  $n\geq
\frac{(p+1)}{2}+1$  if $p$ is odd. Similar, using Lemma \ref{Jost-asymptotics}, we get the leading orders of the resonances. In the even case the resonances can converge to complex number - zeros of the factor $(\l-\tilde{v}_{p-1})\left\{(\l-\tilde{v}_p)-
\frac{2\delta}{\l-v}\right\}-1$ in (\ref{complex-resonances1}) or equivalently zeros of the polynomial $(\l-\tilde{v}_{p-1})(1-q_p(\l-v))-(\l-v)$ (see
(\ref{complex-resonances2})).

This  implies Theorem \ref{th-a0}.

Theorem \ref{th_largeq} follows from \ref{Jost-asymptotics} and iv) in Lemma \ref{l-asymptotics-F}.
\\ \\

In the next Lemma we state the crucial properties of the function
$F.$

\begin{lemma}\label{generalresults}
i) Suppose that $\l_1\in\gamma_k^+,$ for $k=0,$ $1$ or $2,$ and
either\\
a) $f_0^+(\l_1)=0,$ i.e.
 $\l_1$ is an eigenvalue of $J$ with the eigenfunction
 $y_n=f_n^+(\l_1),$ or\\
 b) $\l_1=v.$ Let $\l_1$ also denote the projection of
 $\l_1\in\l_k^+$ on $\C.$

 Then $(-1)^k\dot{F}(\l_1) <0$ and function
$F$ has simple zeros at all bound states of $J.$ Moreover if
$\l_1=v,$  then $\tilde{\vp}_0(v)=0,$ $f_0^+(\l)$ is analytic at
$\l=v\in\gamma_1^+$ and $f_0^+\neq 0.$

ii)  We have $F(\l)=$
\[\lb{kvform}
\varphi_2\left(\tilde{\vartheta}_0+
\frac{\phi}{\varphi_2}\tilde{\varphi}_0\right)^2
+\frac{1-\Delta^2}{\varphi_2}\tilde{\varphi}_0^2=\varphi_2\left(\tilde{\vartheta}_0+
\frac{\phi}{\varphi_2}\tilde{\varphi}_0\right)^2
+\frac{-(\l^2-\l_0^2)(\l^2-\l_1^2)}{4a^2\varphi_2}\tilde{\varphi}_0^2,\]
where $\l_0=\l_0^\mp$ and $\l_1=\l_1^\pm$ are the endpoints of
$\sigma_{\rm ac}(J^0)=[\l_{0}^+,\l_1^-]\cup [\l_{1}^+,\l_0^-].$

We have $F(\l) <0,$ for $\l\in (\l_{0}^+,\l_1^-),$ and  $F(\l) >0,$
for $\l\in (\l_{1}^+,\l_0^-).$

\end{lemma}

\begin{figure}[h]
\begin{center}
\includegraphics*[scale=0.90]{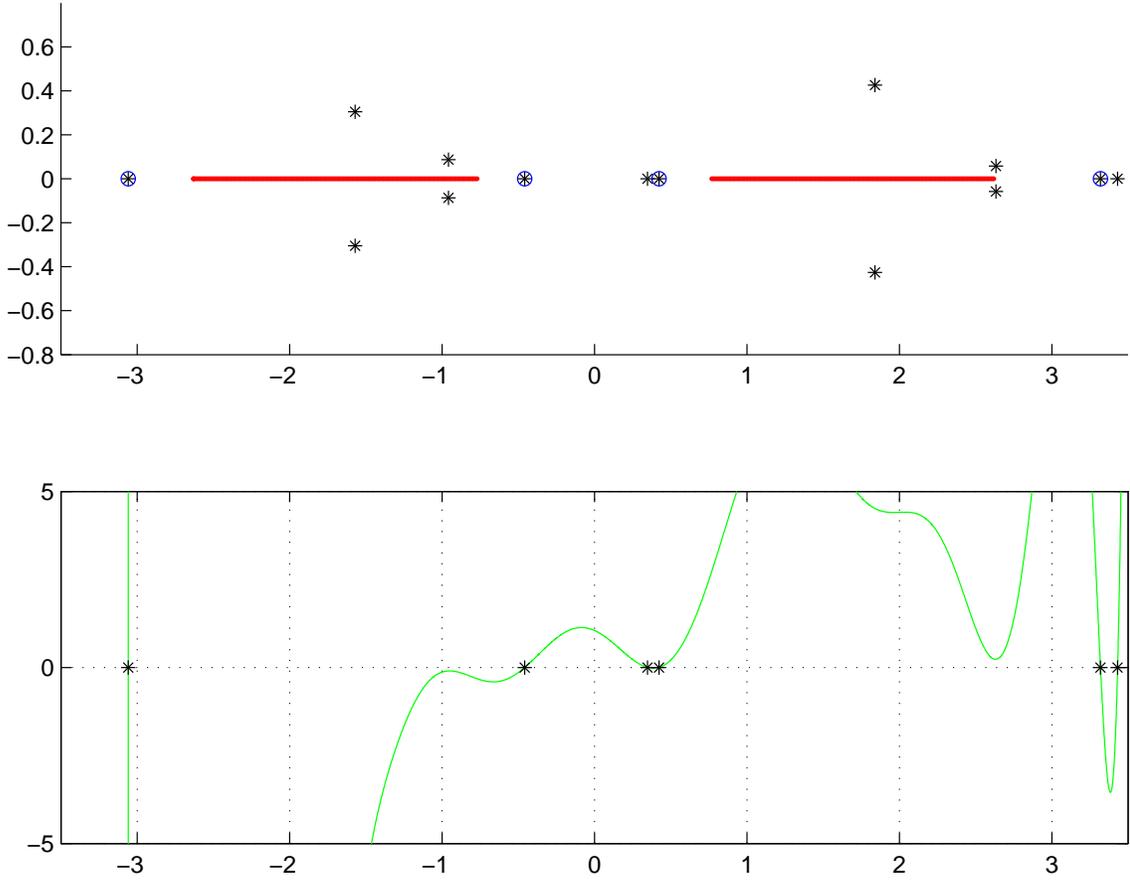}
\end{center}
 \caption{\small Function $F$ and the states, the bound states are encircled. \label{Fig1} }
\end{figure}
{\bf Remarks.} 1) Lemma \ref{generalresults} (with proper modifications) is also
true for general Jacobi operators on the half-line and is proven in
paper \cite{IK3}. The methods remind the approach of \cite{K2} to  the periodic Schr{\"o}dinger
operator plus compactly supported potentials on the half-line.

2) It follows that $F(\l),$ which is real on the real axis, is
decreasing function at any eigenvalue $\l_1\in\gamma_{0,2}^+,$  and
increasing function at any eigenvalue $\l_1\in\gamma_1^+.$

It follows that all bound states of $J$ are simple and that for any
two eigenvalues $\l_{1,2}\in\gamma_k^+$  such that the interval
$\Omega^+=(\l_1,\l_2)\subset\gamma_k^+$ does not contain any other
bound states there is an odd number of antibound states $\l_0$ in
the same interval $\Omega^-\subset\gamma_k^-\subset\Lambda_-$ on the
second sheet, and  $(-1)^k\dot{F}(\l_0)>0.$



3) From ii) it follows that  there is always at least one eigenvalue
in the middle gap $\gamma_1$ or a virtual state at $\l_1^\pm.$
Moreover, using that from Lemma \ref{l-asymptotics-F} it follows that function $F$ has the same sign  when $\l\rightarrow\pm\infty,$ we get that\\
if $\sign(F)(\pm\infty)<0,$  then there are at least two
eigenvalues: one in $\gamma_1$ another in $\gamma_2$ (which can
become virtual states);
 if $\sign(F)(\pm\infty)>0$ and  $0<a<1,$ then there are at least two eigenvalues: one
 in
 $\gamma_0,$ another in $\gamma_1$ (which can
become virtual states).

Now the proof  of Theorem \ref{th-charact-states} follows from the
following two Lemmata which are proved in \cite{IK3} in the general
case.

\begin{lemma}[Virtual states]\lb{L-virt} Let $\l_0$ denote any of $\l_{0,1}^\pm$  and
let $\lambda=\l_0+\epsilon$ for
$\epsilon>0$ small enough.\\
i) If $\l,\l_0\neq v$ and  $f_0^+(\l_0)=0,$ then $\l_0$ is a simple
zero of $F,$  $\l_0$ is a virtual state of $J,$ and
\begin{equation}\lb{3.6}
f_0^+(\l)=\tilde{\varphi}_0(\l_0)c\sqrt{\epsilon} +{\mathcal
O}(\epsilon),\,\,\mR_n(\l)=\frac{f_n^+(\l)}{\tilde{\varphi}_0(\l_0)c\sqrt{\epsilon}}(1+{\mathcal
O}(\sqrt{\epsilon})),\,\,c\tilde{\varphi}_0(\l_0)\neq 0.
\end{equation}
ii) If $\l_0=v$ (which happens if $a=1$) and
$\tilde{\varphi}_0(\l_0)\neq 0,$ then $F(\l_0)\neq 0$ and each
$\mR_n(.),$ $n\in\N,$ does not have singularity at $\l_0$ and
$\l_0$ is
not a virtual state of $J.$\\
iii) If $\l_0=v$  and $\tilde{\varphi}_0(\l_0)= 0,$ then $\l_0$ is
a virtual state of $J,$ $f_0^\pm (\l_0)\neq 0,$  $\l_0$ is simple zero
of $F,$ and each $(\mR_n(.))^2,$ $n\in\N,$ has pole at $\l_0.$
\end{lemma}

\begin{lemma} \lb{L-ZerosF} The projection $\pi:\,\,\Lambda \mapsto \C$ of the set of states of  $J$ on $\Lambda$ coincides with the set of zeros
of $F$ on the complex plane $\C:$
$$\pi\gS\,(J)={\rm Zeros}\,(F).$$
Moreover, the multiplicities of bound states and resonances are
equal to the multiplicities of zeros of $F.$ All bound states are
simple.

Suppose  $\l_0=\l_0^\pm$ or $\l_0=\l_1^\pm$ and
$\tilde{\varphi}_0(\l_0)\neq 0.$ Then $\l_0$ is a virtual state iff $F(\l)$ has  zero at $\l_0.$ It will be automatically
simple.
\end{lemma}

In the next Lemma we consider a special case when we have a simple criterium when $\l=v$ is a
state.

 \begin{lemma}\label{Prop-l=v}

  Suppose that $p$ is even and for any $n\in\N,$
 $\tilde{v}_{2n+1}=v_{2n+1}=v.$

  Then $\tilde{\varphi}_0(v)=0$ and $F(v)=0.$ Thus $\l=v$ is a
  state.

\end{lemma}
{\bf Proof:} From the well known explicit formula
$\varphi_{2n}=(\l-v)\sin n2\varkappa /\sin 2\varkappa$ it follows
that $\varphi_{2n}(v)=0$ for any $n\in\N.$ From the equation $Jy=\l y$ we
have the iteration formula:
\begin{equation}\lb{it}\tilde{\varphi}_{2n}=
\frac{(\l-\tilde{v}_{2n+1})\tilde{\varphi}_{2n+1}-\tilde{\varphi}_{2n+2}}{a}.
\end{equation}
But $\tilde{\varphi}_{k}=\varphi_{k}$ for $k\geq p+1$ and
$\tilde{v}_{2n+1}=v.$ Thus, starting with $2n=p$ and iterating
(\ref{it}), we get $\tilde{\varphi}_{2n}(v)=0,$
$\tilde{\varphi}_{2n-2}(v)=0,$ $\ldots,$ $\tilde{\varphi}_{0}(v)=0:$
all functions $\tilde{\varphi}_k$ with even indexes are zeros at
$\l=v.$ Then from (\ref{formulaF}) it follows that $F(v)=0.$ \qed

\section{Examples $p=1$ and $p=2.$}\label{S-examples}
\setcounter{equation}{0}
In this Section we consider the
 special cases $p=1$ and $p=2$ when the properties of the states can be analyzed in more details.
Using that
\begin{align*}
f_{p}^+=&{(\l-\tilde{v}_{p+1})f_{p+1}^+-a_{p+1}f_{p+2}^+\/a_p}
=\frac{\x^{p+1}_+}{a_p}((\l +v)-am_+)
=\frac{\x^{p+1}_+}{a_p}\frac{1+a\xi^2_-}{\l -v},\\
f_{p-1}^+&={(\l-\tilde{v}_{p})f_{p}^+-a_{p}f_{p+1}^+\/a_{p-1}}=\frac{\x^{p+1}}{a_pa_{p-1}}
\left( (\l-\tilde{v}_p)\frac{1+a\xi^2_-}{\lambda -v}-a_p^2\right),
\end{align*}
and $\varphi_2=\l -v,$  $\xi_+^2\xi_-^2=1,$
$\xi^2_++\xi_-^2=2\Delta=(\l^2-v^2-a^2-1)/a$ $\Rightarrow$ $2a\Delta
+1+a^2=\l^2 -v^2,$ we get for $p=1,$ $a_0=a,$ $a_1=1,$ $v_1=v+q_1,$
\begin{align*}
F(\l)&=\varphi_2f_0^+f_0^-=a^{-2}\left(-q_1\l^2 +\l \left[
q_1^2+a^2\right]+\left(q_1^2v+q_1(v^2+1-a)-va^2\right)\right),
\end{align*}
If $v=v_1$ (the unperturbed case) then $F=a^2(\l -v)$ and $\l=v$ is the
only state, see Lemma \ref{Prop_states_nonperturbed}. The discriminant of the quadratic equation is
\begin{align*}
D&=(a^2-(v_1^2-v^2))^2+4(v_1-v)^2=(q_1^2 +2vq_1 -a^2)^2
+4q_1^2>0\,\,\mbox{if}\,\, v_1\neq v.
\end{align*}
Thus we get that the states are real. By Lemma \ref{generalresults}
 on $F,$ part ii), both states are bound states: no resonances for
$p=1.$ One can check directly that if perturbation is non-trivial
($v\neq v_1$) then there are no virtual states if $p=1:$
 if $\l_0=\l_0^\pm$ is virtual state then
$\xi^2_-=\Delta(\l_0)=\pm 1,$ and
$$f_0^+(\l_0)=\frac{\pm 1}{a}\left((\l-v_1)\frac{1\pm a}{\l-v}
-1\right)=0\,\,\Leftrightarrow\,\,(\l_0-v_1)(1\pm a)=\l_0
-v,\,\,\l_0\neq v,$$ which never happens. Thus we have
\begin{proposition}\label{Prop_p1} For $p=1,$ $v\neq v_1,$ $J$ has two real bound
states:
$$\l_{\pm}=\frac{[q_1^2+a^2]\pm\sqrt{(q_1^2 +2vq_1
-a^2)^2 +4q_1^2}}{2q_1}=\frac{q_1}{2}+\frac{a^2}{2q_1}\pm
\sqrt{\left(\frac{q_1}{2}+v-\frac{a^2}{2q_1}\right)^2+1}.$$
\end{proposition}
In the limit $a\rightarrow 0,$ we get straightforward
$\l_{1,2}=\frac12(q_1\pm\sqrt{(2v+q_1)^2 +4})+{\mathcal O}(a^2).$ As
$v\rightarrow\infty,$ we have $\l_{1,2}\sim \pm
v\rightarrow\pm\infty.$ As $v\rightarrow 0,$ we have
$\l_{1,2}\rightarrow (2q_1)^{-1}([q_1^2+a^2]\pm\sqrt{(q_1^2 -a^2)^2
+4q_1^2}).$ Next we get:\\if $q_1\rightarrow 0,$ then $\l_+\sim
a^2/q_1\rightarrow \infty,$ and $\l_-\rightarrow v;$\\ if
$q_1\rightarrow\infty,$ we have $\l_+\sim q_1\rightarrow\infty,$ and
$\l_-\rightarrow -v.$








%

%

\vspace{4mm}

Now we  consider in detail the properties of the states
in the simplest non-trivial case $p=2,$ which allows the complex
resonances.
Let  $D(p_3)$ denote  the generalized  discriminant of a special cubic polynomial which will be explained below and  given by the following cumbersome
formula:
\begin{align}
D(p_3)=&(vq_2+q_2^2)^2q_2^2(2vq_2-v^2-a^2-1)^2
-4(vq_2+q_2^2)^3\{(vq_2-v^2-1)(vq_2-a^2)-v^2a^2\}-\nonumber\\
&-4q_2^4(2vq_2-v^2-a^2-1)^3+\nonumber\\&+18q_2(vq_2+q_2^2)q_2(2vq_2-v^2-a^2-1)\{(vq_2-v^2-1)(vq_2-a^2)-v^2a^2\}-\nonumber\\
&-27q_2^2\{(vq_2-v^2-1)(vq_2-a^2)-v^2a^2\}^2. \label{cumbersome}
\end{align}

\begin{proposition}\lb{Prop_p=2} i) Suppose $p=2$ and $q_2\neq 0.$ Then $J$ has always two
bound states and two resonances (or virtual states). In the limit
$a\rightarrow 0+,$ the bound states converge to

\begin{equation}\lb{lambda12}
\l_{1,2}^0=\pm\sqrt{\left(v+\frac{q_1-q_2}{2}\right)^2+1}+\frac{q_1+q_2}{2}
\end{equation}
and the resonances converge to \begin{equation}\lb{lambda34}
\l_{3,4}^0=v+\frac{q_1}{2}\pm\sqrt{\frac{q_1(q_2q_1-4)}{4q_2}}.
\end{equation}
 Suppose
that $\tilde{v}_1=v.$ Then $\l=v$ is always a state.

 Moreover, let
$D(p_3)$ denote the generalized discriminant given by Formula
(\ref{cumbersome}). Then, all  four states of $J$ are real iff $D(p_3)>0.$ If $D(p_3)<0,$ then there are always two complex
conjugated resonances.

ii) Suppose $p=2,$ $\tilde{v}_1=v.$ We have the following asymptotic
properties of the states:\\
   1) for $q_2$ small
enough, $J$ has precisely two non-real complex
conjugated resonances;\\
 2) in the limit $q_2\rightarrow \infty,$  the states
 of $J$ either go to infinity or converge to the real state $\l =v;$\\
 3) in the limit $v\rightarrow\infty,$ the states are of order $|v|.$
 Moreover, let $\mu_{1,2,3}$ denote the zeros of $\mu^3-\mu^2-\mu+q_2=0,$
 which  are real if
 \begin{equation}\lb{if_real}
 q_2\in \left[ \frac{11-\sqrt{11^2+5\cdot 27}}{27},
 \frac{11+\sqrt{11^2+5\cdot 27}}{27}\right]
 \end{equation}
 and contain one complex
 conjugate pair otherwise. Then $\l_{1,2,3}/v\rightarrow
 \mu_{1,2,3}$ in the limit $v\rightarrow\infty;$
 4) for $v$ small enough we do not have non-numerical results;\\
 5) in the limit $a\rightarrow 0+,$ the two resonances converge to
 $\l=v.$
\end{proposition}
Some of this results can be generalized to any $p$ (see Theorem
\ref{th_largeq}). The fact that $\l=v$ is always a state in the
special case $p=2,$ $\tilde{v}_1=v,$ can be generalized to any even
$p$ (see Lemma \ref{Prop-l=v}).
We proceed now to the proofs.\\
{\bf Proof:}
From the properties of function $F$ we know that $J$ has always at
least two bound states (or eventually virtual states).

For $p=2$ we have $$f_0^+(\l)={\x^2_+\/a^2}\rt(
(\l-v_{1})((\l-v_2)a-a^2m_+)-a\rt),\,\,f_0^-(\l)=\overline{f_0^+(\overline{\l})}.$$

We get $F=\varphi_2f_0^+f_0^-=$
\begin{align}
&=\frac{1}{a^2}\left[(\l-v)\left\{(\l-v_1)(\l-v_2)-1\right\}^2-\left\{(\l-v_1)(\l-v_2)-1\right\}(\l-v_1)(\l^2-v^2+a^2-1)+\right.\nonumber\\
&\hspace{11cm}+\left.(\l-v_1)^2 a^2(\l+v)\right],\label{this}
\end{align}
where we used that $m_++m_-=2\Phi/\varphi_2$ and
$m_+m_-=-\vartheta_3/\varphi_2,$ where $\varphi_2=\l-v,$
$\Phi=(\l^2-v^2+a^2-1)/2a,$ $\vartheta_3=-\l -v.$

Suppose $v_1=v$ and $\l\neq v.$ From (\ref{this})  we get
${\displaystyle \frac{a^2 F(\l)}{\l-v}=}$
\begin{align}
&=-(\l^2-v^2-a^2-1)(\l-v)q_2+(\l-v)^2q_2^2+a^2=\label{detr}\\
&=-q_2\l^3+\l^2(vq_2+q_2^2)-\l
q_2(2vq_2-v^2-a^2-1)+(vq_2-v^2-1)(vq_2-a^2)-v^2a^2 \lb{cubic}
\end{align}

 As $F(\l)=\frac{\l
-v}{a^2}p_3(\l)$ we have: $\sign F=-\sign q_2$ in the limit $\l\rightarrow\pm\infty.$ As $F$
is strictly negative under the first band and strictly
positive over the second band, we have:\\
if $q_2<0$ then there is one bound state in $\gamma_0^+;$
if $q_2>0$ then there is one bound state in $\gamma_2^+;$\\
state $\l=v\in\gamma_1$;\\
the other 2 states are either real, then they belong to the same
gap, or complex conjugate.

 The right hand side of (\ref{cubic}) is the cubic polynomial with real coefficients in
$\l$ and can have $3$ real zeros or one real zero and two complex
conjugated zeros. Denote the respective coefficients in
(\ref{cubic}) by $k_0,k_1,k_2,k_3,$ then we have
\begin{align*}
&\frac{a^2 F(\l)}{\l-v}=p_3(\l)=k_0\l^3+k_1\l^2+ k_2\l+k_3,\\
& k_0=-q_2,\,\,k_1=vq_2+q_2^2,\,\,k_2=-q_2(2vq_2-v^2-a^2-1),\,\,
k_3=(vq_2-v^2-1)(vq_2-a^2)-v^2a^2.
\end{align*}

 Remark 3.104 on page 127 from \cite{Vi}, states that if the
generalized  discriminant of $p_3$
\begin{equation}\label{cubicdiscr} D(p_3)=k_1^2k_2^2
-4k_1^3k_3-4k_0k_2^3+18k_0k_1k_2k_3-27k_0^2k_3^2\end{equation} is
strictly positive then all the zeros of (\ref{cubic}) are real and
disctinct. If $D(p_3)<0,$ then there are two complex conjugated
zeros. The discriminant $D(p_3)$ is given in (\ref{cumbersome}).

Thus we have proven the first part i) of Proposition \ref{Prop_p=2}
except asymptotics (\ref{lambda12}) and (\ref{lambda34} which we
postpone to the end of this section.


Suppose $q_2\rightarrow 0.$ Then
$$ k_0=-q_2,\,\,k_1=vq_2+{\mathcal O}(q_2^2),\,\,k_2=q_2(v^2+a^2+1)+{\mathcal O}(q_2^2),\,\,
k_3=a^2+{\mathcal O}(q_2).$$

 Then $D=-27k_0^2k_3^2+{\mathcal O}(q_2^3)=-27q_2^2a^4+{\mathcal O}(q_2^3)$ which implies that there are two non-real resonances.

Suppose $q_2\rightarrow\infty.$ Then directly from the equation
$$-q_2\l^3 +q_2^2\l^2-2vq_2^2\l +v^2q_2^2={\mathcal
O}(q_2)\,\, \Leftrightarrow\,\,\l^2-2v\l +v^2=(\l-v)^2={\mathcal
O}(\l^3q_2^{-1}),\,\,q_2\rightarrow\infty,$$ we get that the states
which remain bounded as $q_2\rightarrow\infty$ converge to $ v.$

Suppose $v\rightarrow\infty.$ Then the equation is $-q_2\l^3
+vq_2\l^2+q_2v^2\l -v^3q_2^2 ={\mathcal O}(v^2).$ Put $\l=v\mu,$
then $-q_2v^3(\mu^3 -\mu_2-\mu+q_2)={\mathcal
O}(v^2)\,\,\Leftrightarrow\,\,\mu^3 -\mu^2-\mu+q_2={\mathcal
O}((vq_2)^{-1}),\,\,v\rightarrow\infty.$ The equation $\mu^3
-\mu^2-\mu+q_2=0$ has the generalized discriminant
$D=1+4q_2+4+16q_2-27q_2^2=-27q_2^2+22q_2+5,$ whose
 zeros are
$$x_\pm =\frac{11\pm\sqrt{11^2+5\cdot 27}}{27}.$$
Denote $\mu_{1,2,3}$ the zeros of the equation $\mu^3
-\mu^2-\mu+q_2=0.$ We proved that, as $v\rightarrow\infty,$  the
states $\l_{1,2,3}={\mathcal O}(v),$ moreover  $\l_{1,2,3}/v$
converge  to the zeros of the equation $\mu^3 -\mu^2-\mu+q_2=0,$
which are real if $q_2\in [x_-,x_+].$ If $q_2<x_-$ or $q_2 >x_+$
then the two zeros are complex conjugated. Thus we have proven the
part ii) of Proposition
 \ref{Prop_p=2}.
In the case $v\rightarrow 0,$ the equation does not simplify.

{\bf Asymptotics (\ref{lambda12}) and (\ref{lambda34}) as
$a\rightarrow 0.$} As a special case of the asymptotics in Lemma
\ref{Jost-asymptotics}, we get, for $\l\in\gamma_1^+,$
\begin{align*}
f_0^+&=\frac{1}{a^2}\left(\frac{a}{2\delta} +{\mathcal
O}(a^3)\right)\left((\l-v_1)\left[(\l-v_2)a +{\mathcal
O}(a^3)\right]-a\right)=\\
&=\frac{1}{a^2}\left(\frac{(\l-v_1)(\l-v_2)-1}{2\delta}a^2+{\mathcal
O}(a^4)\right)=\frac{(\l-v_1)(\l-v_2)-1}{\l^2-v^2-1}+{\mathcal
O}(a^2).
\end{align*}
Thus the bound states in $\gamma_1^+$ which are solutions  of
$f_0^+(\l)=0$ in the limit $a\rightarrow 0$ are asymptotically
solutions of the  equation
$$\frac{(\l-v_1)(\l-v_2)-1}{\l^2-v^2-1}=0\,\,\Leftrightarrow\,\,\l^2-(v_1+v_2)\l +v_1v_2-1=0$$ if $\l^2\neq v^2+1$
(which happens if $\l$ is a virtual state). We get two solutions
($\tilde{v}_1=v+q_1,$ $\tilde{v}_2=-v+q_2$)
$$z_{1,\pm}=\frac{\tilde{v}_1+\tilde{v}_2\pm\sqrt{(\tilde{v}_1-\tilde{v}_2)^2+4}}{2}=\pm\sqrt{\left(v+\frac{q_1-q_2}{2}\right)^2+1}\pm\frac{q_1+q_2}{2}$$ which are the
leading terms in the expansion of the bound states in $\gamma_1^+$ as
$a\rightarrow 0$ (see Theorem \ref{th-a0}).

Similarly we get the resonances in $\gamma_1^-$ and states in
$\gamma_{0,2}^\pm$ in the limit $a\rightarrow 0.$ The resonances in $\gamma_1^-$ are formally also zeros of
$f_0^-(\l)$ in $\gamma_1^+$ which, if $v_1\neq v,$ in the leading order $a^{-2}$ are solutions of the equation
$
(v_2+v)\l^2-(v_2v+v^2+v_1(v_2+v))\l
+v_1(v_2v+v^2+1)-v=0$ with zeros given in  (\ref{complex-resonances0}):
$$
\mu_{1,2}^0
=\frac{q_2(2v+q_1)\pm\sqrt{q_2q_1(q_2q_1-4)}}{2q_2}=v+\frac{q_1}{2}\pm\sqrt{\frac{q_1(q_2q_1-4)}{4q_2}}
$$
 which can be real antibound states or  complex conjugated
resonances. If  $v_1=v$ then
the leading order as $a\rightarrow 0$ of the antibound state is
$\l=v.$

\qed

Suppose now that $\l_0$ is a real double root:
$F(\l_0)=0,\,\,\dot{F}(\l_0)=0.$ Suppose $\l_0\neq v.$ Using the
identity
$$a^2\dot{F}(\l)=a^2\frac{F(\l)}{\l -v}
+(\l-v)\frac{\partial}{\partial\l}\left(a^2\frac{F(\l)}{\l-v}\right),$$
we have at $\l=\l_0:$
$$a^2\dot{F}(\l_0)=(\l_0-v)\frac{\partial}{\partial\l}\left(a^2\frac{F(\l)}{\l-v}\right)_{|\l=\l_0}=0.$$

In the special case $p=2,v_1=v,$ we get using (\ref{cubic}):
\begin{equation*}\lb{cubicderivative}
\frac{\partial}{\partial\l}\left(a^2\frac{F(\l)}{\l-v}\right)_{|\l=\l_0}=
-3q_2\l_0^2+2q_2(v+q_2)\l_0- q_2(2vq_2-v^2-a^2-1)=0.
\end{equation*}
It follows that if $F(\l_0)=\dot{F}(\l_0)=0$ and $\l_0\neq v$ then
$\l_0$ is a zero of the quadratic equation $3\l^2-2(v+q_2)\l+
(2vq_2-v^2-a^2-1)=0:$
$$\l_0=\frac{(v+q_2)\pm\sqrt{(v+q_2)^2-3(2vq_2-v^2-a^2-1)}}{3}=\frac{(v+q_2)\pm\sqrt{(2v-q_2)^2+3(a^2+1)}}{3}.$$ The state
 has multiplicity $2$ and necessarily are the antibound states as all
bound states  are simple. Thus we get
\begin{proposition}[Double antibound state] Suppose $p=2,$ $v_1=v.$ Operator $J$ has
precisely one antibound state at $\l_0$ of multiplicity $2$  iff the discriminant given in (\ref{cumbersome}) is zero.
Moreover, the double antibound state different from $\l=v$ is given
by either of the two following formul{\ae}:
$$\frac{(v+q_2)\pm\sqrt{(2v-q_2)^2+3(a^2+1)}}{3}.$$ The state
$\lambda=v$ is always simple if $a\neq 0.$

\end{proposition}
\vspace{1cm}
 {\bf Acknowledgement.} The authors thank the referee for remarks.

\end{document}